\documentclass[10pt,twocolumn,draft]{article} 
\usepackage{latex8}
\usepackage{latexsym,amssymb,amsmath}
\usepackage{times,mathptm}
\usepackage[thmmarks]{ntheorem}
\usepackage[latin1]{inputenc}

\def\mms#1{{\imm{
\mathord{
\mathchoice{\hbox{#1}}{\hbox{#1}}{\hbox{\tiny #1}}{\hbox{\tiny #1}}
}}}}

\theoremnumbering{arabic}
\theoremstyle{plain}
\theoremsymbol{}
\makeatletter
\theorembodyfont{\itshape}
\theoremheaderfont{\normalfont\bfseries}
\theoremseparator{}
\newtheorem{thm}{\quad Theorem}[section]
\newtheorem{lem}[thm]{\quad Lemma}
\newtheorem{prop}[thm]{\quad Proposition}

\newtheorem{cor}[thm]{\quad Corollary}
\theorembodyfont{\upshape}
\newtheorem{defn}[thm]{\quad Definition}
\newtheorem{ex}[thm]{\quad Example}
\newtheorem{rem}[thm]{\quad \em Remark}
\theoremstyle{nonumberplain}
\theoremheaderfont{\itshape}
\theorembodyfont{\normalfont}
\theoremsymbol{\ensuremath{\Box}}
\newtheorem{proof}{\quad Proof.}
\qedsymbol{\ensuremath{\Box}}

\leftmargin\leftmargini
\labelwidth         \z@
\topsep \z@

\def\@listI{\leftmargin\leftmargini} \@listI
\def\@listi{\leftmargin\leftmargini  \topsep \z@ plus 1pt minus 1pt}
\def\@listii{\leftmargin\leftmarginii\labelwidth\leftmarginii
    \advance\labelwidth-\labelsep    \topsep \z@}

\def\section{\@startsection {section}{1}{\z@}
   {8pt plus 2pt minus 2pt}{8pt plus 2pt minus 2pt} {\large\bf}} 
\def\subsection{\@startsection {subsection}{2}{\z@}
   {11pt plus 2pt minus 2pt}{11pt plus 2pt minus 2pt} {\elvbf}}

\makeatother

\let\imm\ensuremath
\let\nmodels\nvDash
\let\phi\varphi

\def\boeta{{{\bf\eta}}}
\def\one{{{\bf 1}}}

\def\s{\mms{\bf s}}
\def\N{{\mathbb N}}

\def\R{{\mathbb R}}
\def\up{{\uparrow}}
\def\dn{{\downarrow}}

\let\impl\to

\let\mmodels\models
\def\G{{\mathbf{G}}}
\def\LC{{\bf LC}}
\def\IPL{{\bf IPL}}

\def\Frm{{\rm Frm}}

\def\Distr{{\rm Distr}}
\def\I{{\imm{\mathfrak{I}}}} \def\IL{{\bf IL}} 
\def\J{{\imm{\mathfrak{J}}}} 
 
\def\Gd{{\G_V}}
\def\Q{{\ensuremath{\mathord{\mathsf{Q}}}}}
\def\T{{\bf\sf T}}
\def\lo{\mathbf{L}}
\newcommand\tup[1]{\overline{#1}}
\newcommand{\modelsg}{\mathrel{\Gd \models}}
\newcommand{\modelsl}{\mathrel{\lo \models}}
\newcommand{\IMPLIES}{\ \Rightarrow\ }

\def\bbR{{\mathbb{R}}}
\def\bbQ{{\mathbb{Q}}}
\def\succ{K}

\title{Characterization of the Axiomatizable\\ Prenex Fragments of
First-Order Gödel Logics}

\author{Matthias Baaz\thanks{Research supported by FWF grant
P15477--MAT}\\ 
Technische Universität Wien\\ 
A--1040 Vienna, Austria\\
baaz@logic.at \and
Norbert Preining$^*$\\
Technische Universität Wien\\ 
A--1040 Vienna, Austria\\
preining@logic.at \and
Richard Zach\\
University of Calgary\\
Calgary, AB T2N 1N4,
Canada\\ rzach@ucalgary.ca}

\begin{document}

\maketitle

\begin{abstract} 
The prenex fragments of first-order infinite-valued Gödel logics are
classified.  It is shown that the prenex Gödel logics characterized
by finite and by uncountable subsets of $[0, 1]$ are axiomatizable,
and that the prenex fragments of all countably infinite Gödel
logics are not axiomatizable.
\end{abstract}

%%%%%%%%%%%%%%%%%%%%%%%%%%%%%%%%%%%%%%%%%%%%%%%%%%%%%%%%%%%%%%%%%%%%%%%%%
%
% INTRODUCTION
%
%%%%%%%%%%%%%%%%%%%%%%%%%%%%%%%%%%%%%%%%%%%%%%%%%%%%%%%%%%%%%%%%%%%%%%%%%

\vspace*{-5pt}
\section{Introduction}

Gödel logics are one of the oldest and most interesting families of
many-valued logics.  Introduced by Gödel in \cite{Godel:32}, they
provide the first examples of intermediate logics
(intermediate, that is, in strength between classical and
intuitionistic logics).  Dummett \cite{Dummett:59} was the first to
study infinite valued Gödel logics, axiomatizing the set of
tautologies over infinite truth-value sets by intuitionistic logic
plus the linearity axiom $(A \impl B) \lor (B \impl A)$.  In terms of
Kripke semantics, the characteristic linearity axiom picks out those
accessibility relations which are linear orders.

Gödel logics have recently received increasing attention, both in
terms of foundational investigations and in terms of applications.  One
of the most surprising recent results is that whereas there is only
one infinite-valued propositional Gödel logic, there are infinitely
many different logics at the first-order and already when only ``fuzzy
quantifiers'' are added to the language
\cite{BaazLeitschZach:96,BaazVeith:97}.  In light of the fact that
first-order infinite-valued \L ukasiewicz logic is not axiomatizable,
it is perhaps also surprising that at least one infinite-valued Gödel
logic \emph{is} r.e. \cite{horn,TT}.

Our aim in the present paper is to characterize the axiomatizable
first-order \emph{prenex} Gödel logics, i.e., those truth-value sets
whose first-order validities in prenex form are r.e. This is a first
step toward the characterization (in terms of axiomatizability) of
first-order Gödel logics in general.  Our result is that there is only
one axiomatizable infinite-valued first-order prenex Gödel logic; it
is characterized by any closed uncountable subset of $[0, 1]$.  In
fact, we give an axiomatization based on a version of Herbrand's
theorem for such truth-value sets, and then show that any countably
infinite truth-value set has a set of prenex validities which is not
r.e.

%%%%%%%%%%%%%%%%%%%%%%%%%%%%%%%%%%%%%%%%%%%%%%%%%%%%%%%%%%%%%%%%%%%%%%%%%
%
% SYNTAX AND SEMANTIC
%
%%%%%%%%%%%%%%%%%%%%%%%%%%%%%%%%%%%%%%%%%%%%%%%%%%%%%%%%%%%%%%%%%%%%%%%%%
\section{Syntax and semantics}

First-order Gödel logics are given by truth functions for the
connectives and quantifiers, and a set of truth values.  We work in a
standard first-order language with variables ($x$, $y$, $z$, \dots),
constants ($a$, $b$, $c$, \dots), function symbols ($f$, $g$, \dots),
predicate symbols ($P$, $Q$, $R$, \dots), the predicate constant
$\bot$, connectives ($\land$, $\lor$, $\impl$) and quantifiers
($\exists$, $\forall$). $\lnot A$ is defined as $A\impl\bot$. The set
of well-formed formulas is denoted by~$\Frm(L)$.  The sets of truth
values for the systems we consider are closed subsets of $[0,1]$
containing both $0$ and $1$.  Interpretations are defined as usual:

\begin{defn}
Let $V \subseteq [0,1]$ be some set of truth values which contains
$0$~and~$1$ and is closed in $\bbR$.  A {\em many-valued
interpretation $\I = \langle D, \s\rangle$ based on $V$} is given by a
{\em domain} $D$ and a {\em valuation function} $\s$ which maps
$n$-ary relation symbols to functions $D^n \to V$, $\s(\bot) = 0$,
$n$-ary function symbols to functions from $D^n$ to $D$, and constants
of $L^\I$ and variables to elements of~$D$. $L^\I$ is $L$ extended by
constant symbols for all $d \in D$; if $d \in D$, then $\s(d) = d$.

$\s$ can be extended in the obvious way to a function on all terms in
$L^\I$.  The valuation of formulas in $L^\I$ is defined by:

(1) $\I(P(t_1, \ldots, t_n)) = \s(P)(\s(t_1), \ldots, \s(t_n))$.

(2) $\I(A \land B) =  \min(\I(A), \I(B))$.

(3) $\I(A \lor B) =  \max(\I(A), \I(B))$.

(4) $\I(A \impl B) = \I(B)$ if $\I(A) > \I(B)$ and $= 1$ otherwise.

\noindent Since we defined $\neg A \equiv A \impl \bot$, we also have
$\I(\neg A) = 0$ if $\I(A) > 0$ and $= 1$ otherwise.

The set $\Distr_\I A(x) = \{\I(A(d)) : d \in D \}$ is called the {\em
distribution} of $A(x)$.  The valuations of quantified formulas are
defined by infimum and supremum of their distributions.

(5) $\I(\forall x\, A(x)) = \inf \Distr_\I A(x)$.

(6) $\I(\exists x\, A(x)) = \sup \Distr_\I A(x)$.

\noindent \I{} {\em satisfies} a formula $A$, $\I \mmodels A$, if
$\I(A) = 1$.
\end{defn}
 
\begin{defn} 
The {\em Gödel logic}~$\G_V$ based on a set of truth values $V$ is the
set of all $A \in \Frm(L)$ s.t.\ $\Gd \models A$, i.e., for every \I{}
based on $V$, $\I \models A$.  The logics $\G_\R$, $\G_\dn$, $\G_\up$,
$\G_m$ are based on the truth value sets
\[\begin{array}{ll}
V_\R = [0,1] & 
V_\dn = \{\frac{1}{k} : k \ge 1\} \cup \{0\} \\ 
V_\up = \{1 - \frac{1}{k} : k \ge 1\} \cup \{1\} & 
V_m =  \{0, \frac{1}{2}, \frac{2}{3}, \ldots, \frac{m-2}{m-1}, 1\}.
\end{array}\]
\end{defn}

%%%%%%%%%%%%%%%%%%%%%%%%%%%%%%%%%%%%%%%%%%%%%%%%%%%%%%%%%%%%%%%%%%%%%%%%%
%
% RELATION BETWEEN GÖDEL LOGICS
%
%%%%%%%%%%%%%%%%%%%%%%%%%%%%%%%%%%%%%%%%%%%%%%%%%%%%%%%%%%%%%%%%%%%%%%%%%

\section{Relationships between Gödel logics}

In the propositional case, the relationships between finite and
infinite valued Gödel logics are well understood.  Any choice of an
infinite set of truth-values results in the same propositional Gödel
logic, viz., Dummett's \LC~\cite{Dummett:59}.  Furthermore, we know
that \LC{} is the intersection of all finite-valued propositional
Gödel logics, and that it is axiomatized by intuitionistic
propositional logic~\IPL{} plus the schema $(A \impl B) \lor (B \impl
A)$.  \IPL{} is contained in all Gödel logics (finite- or
infinite-valued).

In the first-order case, the relationships are somewhat more involved.
First of all, let us note that Intuitionistic predicate logic~\IL{} is
contained in all first-order Gödel logics, since the axioms and rules
of \IL{} are sound for the Gödel truth functions.
As a consequence, we will be able to use any
intuitionistically sound rules and intuitionistically valid formulas
when working in any of the Gödel logics.

\begin{prop}\label{basic-cont}
(1)~$\G_m \supsetneq \G_{m+1}$,
(2)~$\G_m \supsetneq \G_\up \supsetneq \G_\R$,
(3)~$\G_m \supsetneq \G_\dn \supsetneq \G_\R$.
\end{prop}

\begin{proof}
The only nontrivial part is proving that the containments are strict.
For this note that
\(
(A_1 \impl A_2) \lor \ldots \lor (A_m \impl A_{m+1})
\)
is valid in $\G_m$ but not in $\G_{m+1}$. Furthermore, let
%\begin{eqnarray*}
%C_\up & = & \exists x(A(x) \impl \forall y\, A(y)) {\rm\ and\ } \\
%C_\dn & = & \exists x(\exists y\, A(y) \impl A(x)).
%\end{eqnarray*}
\(
C_\up = \exists x(A(x) \impl \forall y\, A(y)) {\rm\ and\ } 
C_\dn = \exists x(\exists y\, A(y) \impl A(x)).
\)
$C_\dn$ is valid in all $\G_m$ and in $\G_\up$ and $\G_\dn$; $C_\up$
is valid in all $\G_m$ and in $\G_\up$, but not in $\G_\dn$; neither
is valid in $\G_\R$ \cite[Corollary~2.9]{BaazLeitschZach:96}.
\end{proof}

The formulas $C_\up$ and $C_\dn$ are of some importance in the study
of first-order infinite-valued Gödel logics.  $C_\up$ expresses the
fact that every infimum in the set of truth values is a minimum, and
$C_\dn$ states that every supremum (except possibly 1) is a maximum.
The only three quantifier shifting rules which are not
intuitionistically valid are:
%\[
%\begin{array}{@{}cc@{}}
%(\forall x\, A(x) \land B) \equiv \forall x(A(x) \land B) & 
%(B \impl \forall x\, A(x)) \equiv \forall x(B \impl A(x)) \\
%(\exists x\, A(x) \land B) \equiv \exists x(A(x) \land B) &
%(B \impl \exists x\, A(x)) \rimpl \exists x(B \impl A(x)) \\
%(\forall x\, A(x) \lor B) \impl \forall x(A(x) \lor B) &
%(\forall x\, A(x) \impl B) \rimpl \exists x(A(x) \impl B) \\
%(\exists x\, A(x) \lor B) \equiv \exists x(A(x) \lor B) &
%(\exists x\, A(x) \impl B) \equiv \forall x(A(x) \impl B)
%\end{array}
%\]
$$
\begin{array}{rcl}
\forall x(A(x) \lor B)  & \impl & (\forall x\, A(x) \lor B) \\
(B \impl \exists x\, A(x)) & \impl & \exists x(B \impl A(x)) \\ 
(\forall x\, A(x) \impl B) & \impl & \exists x(A(x) \impl B)
\end{array}
\eqno{\begin{array}{r}(S_1)\\(S_2)\\(S_3)\end{array}}
$$ ($x$ is not free in $B$.) Of these, $S_1$ is valid in any Gödel
logic.  $S_2$ and $S_3$ imply $C_\dn$ and $C_\up$, respectively (take
$\exists y\, A(y)$ and $\forall y\, A(y)$, respectively, for $B$).
$S_2$ and $S_3$ are, respectively, both valid in $\G_\up$, invalid and
valid in $\G_\dn$, and both invalid in $\G_\R$.  $\G_\up$ is thus the
only Gödel logic where every formula is equivalent to a prenex
formula.  This also implies that $\G_\up \neq \G_\dn$.  In fact, we
have $\G_\dn \subsetneq \G_\up$; this follows from the following
theorem.

\begin{thm}\label{herbrand}
$\G_\up = \bigcap_{m \ge 2} \G_m$
\end{thm}

\begin{proof}
By Proposition~\ref{basic-cont}, $\G_\up \subseteq \bigcap_{m \ge 2}
\G_m$.  We now show the reverse inclusion.  Since all quantifier
shifting rules are valid in $\G_\up$, any formula~$A$ is equivalent to
a prenex formula.  For any given prenex formula $A \equiv \Q_1
x_1\ldots\Q_n x_n\,B(\bar x)$ ($B$ quantifier free) we may define the
\emph{Herbrand form}~$A^H$ of $A$ as usual as $\exists x_{i_1}\ldots
\exists x_{i_m}\, B(t_1, \ldots, t_n)$, where $\{x_{i_j} : 1 \le j \le
m\}$ is the set of existentially quantified variables in $A$, and
$t_i$ is $x_{i_j}$ if $i = i_j$, or is $f_i(x_{i_1}, \ldots, x_{i_k})$
if $x_i$ is universally quantified and $k = \max \{j : i_j < i\}$. We
will write $B(t_1,\ldots,t_n)$ as $B^F(x_{i_1}, \ldots, x_{i_m})$ if
we want to emphasize the free variables.

\begin{lem}
$\G_\up \models A$ iff $\G_\up \models A^H$.
\end{lem}

\begin{proof}
If: Suppose that $\G_\up \nmodels A$.  In $\G_\up$, every infimum is a
minimum, i.e., if $\inf \Distr_\I A(x) = v$ then for some $d \in D$,
$\I(A(d)) = v$.  Hence, we can extend $\I$ by interpretations for the
$f_i$ as in the classical case.  Only if: Obvious.
\end{proof}

It now suffices to show that if $\G_\up \nmodels A$, where $A$ is existential,
then $\G_m \nmodels A$ for some~$m$.  For this we need the following

\begin{lem}
Given $\I = \langle D, \s\rangle$ and $v<1$, define $\I_v = \langle D,
\s_v\rangle$ where $\s_v(P)(d_1, \ldots, d_n) = \I(P(d_1, \ldots,
d_n))$ if $\I(P(d_1, \ldots, d_n))\le v$ and $= 1$ otherwise.  If $A
\in \Frm(L^\I)$ does not contain $\forall$, then $\I_v(A) = 1$ if
$\I(A) > v$ and $\I_v(A) = \I(A)$ if $\I(A) \le v$.
\end{lem}

\begin{proof}
By induction on the complexity of $A$. \qed\end{proof}

Now suppose there is an interpretation~\I{} s.t. $\I \nmodels A$. Then
clearly $\sup \Distr_\I B(\bar x) = v < 1$ (where $A \equiv \exists
\bar x\, B(\bar x)$) and there are only finitely many truth values
below~$v$ in~$V_\up$, say $v = 1 - 1/k$.  Construct $\I_v$ as in the
previous lemma.  Since $\I(B(\bar d)) \le v < 1$, by the lemma
$\I_v(B(\bar d)) \le v$ and so $\sup \Distr_{\I_v} B(\bar x) \le v$.
But $\I_v$ is a $\G_{k+1}$ interpretation, and $\I_v \nmodels A$.
\end{proof}

\begin{cor}
$\G_m \supsetneq \bigcap_m \G_m = \G_\up \supsetneq \G_\dn \supsetneq
\G_\R$
\end{cor}

One basic but important result is that the set of validities of $\G_V$
only depends on the order type of $V$. Let $V$, $V'$ be two truth
value sets, $X$ a set of atomic formulas, and suppose there is an
order-preserving injection~$f: \{\I(B)\colon B \in X\} \to V'$ which
is so that $f(1) = 1$ and $f(0) = 0$.  (Call any such $f$ a
\emph{truth value injection on $X$}.) By a simple induction on $A$, we
have:

\begin{prop}\label{tviso}
Let $A$ be a quantifier free formula, and $X$ its set of atomic
subformulas.  If $\I$, $\I'$ are interpretations on $V$, $V'$,
respectively, and $f$ is a truth value injection on $X$, then
$f(\I(A)) = \I'(A)$.  
\end{prop}

%%%%%%%%%%%%%%%%%%%%%%%%%%%%%%%%%%%%%%%%%%%%%%%%%%%%%%%%%%%%%%%%%%%%%%%%%
%
% TOPOLOGICAL AND ORDER THEORETIC PRELIMINARIES
%
%%%%%%%%%%%%%%%%%%%%%%%%%%%%%%%%%%%%%%%%%%%%%%%%%%%%%%%%%%%%%%%%%%%%%%%%%

\section{Order theoretical preliminaries}

We will characterize the axiomatizable Gödel logics by topological and
order-theoretic properties of the underlying truth value set. The most
important of these properties as regards axiomatizability is the
existence of a non-trivial dense linear subordering of the truth value
set, i.e., a subset $V' \subseteq V$, $\left| V' \right| \ge 2$,
where for all $u, v \in V'$ with $u < v$ there is a~$w \in V'$ such
that~$u < w < v$. In this section we show that there is such a dense
subordering iff $V$ is uncountable.
%\begin{lem}\label{lm:orderembedding}
%  Let~$\strc AR$ be a countable dense linear ordering without first
%  and last elements. Let~$\strc BS$ be an arbitrary countable linear
%  ordering. Then~$\strc BS$ is isomorphic to a subordering of~$\strc
%  AR$.
%\end{lem}
%\begin{proof} \cite[Theorem~2.5]{rosenstein}. \end{proof}

Dense linear orderings are connected to uncountable sets via perfect
sets. We begin by listing some results about perfect sets
from~\cite{kechris}.  All the following notations, lemmas, theorems
are carried out within the framework of Polish spaces, i.e., separable
completely metrizable topological spaces. For our discussion it is
only necessary to know that any closed subset of~$\bbR$ (and hence,
any truth-value set) is such a Polish space.
\begin{defn}
A \emph{limit point} of a topological space is a point that is not
isolated, i.e., for every open neighbourhood~$U$ of~$x$ there is a
point~$y\in U$ with~$y\neq x$. A space is \emph{perfect} if all its
points are limit points.
\end{defn}

It is obvious that intervals of the real line are perfect, but there
are perfect sets which are not intervals:
\begin{ex}
The set of all numbers in the unit interval which can be expressed in
triadic notation using only 0~and~2 is called the \emph{Cantor set};
it is a perfect subset of~$[0,1]$.
\end{ex}
%A more intuitive way to obtain this set is to start with the unit
%interval, take out the middle third (open interval) and restart this
%process with the lower and the upper third. Repeating this we get
%exactly the Cantor set because the middle third always contains the
%numbers which have a~1 in their triadic expansion. Note that the
%point~1/3 of the Cantor set on~$[0,1]$ can be written in triadic
%notation as $0.02222\ldots$.
\begin{prop}\label{perfect}
For any perfect set $P \subseteq \bbR$ there is a unique partition
of~$\bbR$ into countably many intervals such that the intersection of
$P$ with each interval is either empty, the entire
interval or isomorphic to the Cantor set.
\end{prop}

\begin{proof}
  See \cite{winkler:howmuch}, Proposition~1 and discussion.
\end{proof}

To obtain a connection between uncountable sets and perfect sets we
first note that it is possible to embed the Cauchy space into any
perfect space, which yields
\begin{lem}
  If~$X$ is a nonempty perfect Polish space, then the cardinality
  of~$X$ is~$2^{\aleph_0}$; thus all nonempty perfect subsets of ~$\bbR$
  have cardinality of the continuum.
\end{lem}

\begin{proof}
\cite[Corollary 6.3]{kechris}.
\end{proof}

For the other direction, we want to partition an uncountable set 
into a perfect kernel and a countable rest. This is the well known
Cantor-Bendixon Theorem:  
\begin{thm}[Cantor-Bendixon]
Let~$X$ be a Polish space. Then~$X = P \cup C$, with~$P$ a perfect
subset of~$X$ and~$C$ countable open.  $P$ is called the \emph{perfect
kernel} of~$X$.
\end{thm}

As a corollary we obtain that any uncountable Polish space contains a
perfect set, and therefore has cardinality~$2^{\aleph_0}$. Now we can
prove the central theorem:
\begin{thm}\label{bendixon}
  A truth value set (i.e., a closed subset of~$[0,1]$) is uncountable
  iff it contains a non-trivial dense linear subordering. 
\end{thm}
\begin{proof}
  If: Every countable non-trivial dense linear order has order type
  $\boeta$, $\one+\boeta$, $\boeta+\one$, or $\one+\boeta+\one$
  \cite[Corollary~2.9]{rosenstein}, where $\boeta$ is the order type
  of $\bbQ$.  The completion of any ordering of order type~$\boeta$
  has order type~$\lambda$, the order type of $\bbR$
  \cite[Theorem~2.30]{rosenstein}, thus the truth value set must be
  uncountable.

  Only if: We define a dense linear subordering for any uncountable
  set. In fact we will give a dense linear subordering of the perfect
  kernel.

  Since every perfect subset of the real line is a union of intervals
  and sets isomorphic to the Cantor set (Proposition~\ref{perfect}),
  it suffices to show the claim for those sets. For intervals the
  claim is trivial. Now consider the border points in a Cantor set,
  i.e., points which can only be approximated within the Cantor set
  from above or below but not both.  In the ternary notation these are
  the points with a finite number of~0 or a finite number of~2, i.e.,
  their ternary expansions are either $a=0.a_1a_2\ldots a_n$ or
  $b=0.b_1b_2\ldots b_n2222\ldots$ Each border point can be
  approximated by a sequence of inner points~$a^k$. For the $k$-th
  sequence element approximating a border point~$a$ we get $a^k$
  by appending $2k$~zeros and then a sequence of~$020202\ldots$ at the end
  ($a^k = 0.a_1\ldots a_n(00)^k\overline{02}$).  For the $k$-th
  sequence element approximating a border point~$b$ we define an
  approximating sequence $b^k$ by replacing the ternary expansion
  starting from the~$2k$-th 2 with a sequence of 02's ($b^k =
  0.b_1\ldots b_n(22)^k\overline{02}$).  The set of approximations of
  all border points is a dense subset: If $a^k = 0.a_1\ldots
  a_n(00)^{k}\overline{02}$ and $a^{k+1} = 0.a_1\ldots
  a_n(00)^{k}00\overline{02}$ are adjacent points in the sequence,
  then $a' = 0.a_1\ldots a_n(00)^{k}00022222\ldots$ is a border point
  with $a^{k+1} < a' < a^k$, hence there are infinitely many points
  $a'^\ell$ between~$a^k$ and~$a^{k+1}$ in the subset. Similarly for
  adjacent elements of a $b$-sequence. The set of border points is
  countable, therefore the set containing all the approximation
  sequences is countable and has all the necessary properties.
\end{proof}

Note that for example 1/3 and 2/3 would not be in the dense linear
subordering, because between them there is no point of the perfect
set. We would replace 1/3 by a sequence of inner points approximating
1/3 from below and replace 2/3 by a sequence of inner points
approximating 2/3 from above.

%%%%%%%%%%%%%%%%%%%%%%%%%%%%%%%%%%%%%%%%%%%%%%%%%%%%%%%%%%%%%%%%%%%%%%%%%
%
% AXIOMATIZABILITY RESULTS
%
%%%%%%%%%%%%%%%%%%%%%%%%%%%%%%%%%%%%%%%%%%%%%%%%%%%%%%%%%%%%%%%%%%%%%%%%%

\section{Axiomatizability results}

Throughout this section, $V$ is a truth value set which is either
finite or uncountable. Let $\Gd$ be a Gödel logic with such a truth
value set. We show how to effectively associate with each prenex
formula $A$ a quantifier-free formula $A^\ast$ which is valid in $\Gd$
if and only if $A$ is valid.  The axiomatizability of $\Gd$ then
follows from the axiomatizability of \LC{} (in the infinite-valued
case) and propositional $\G_m$ (in the finite-valued case).  Recall
that $A^H$ stands for the Herbrand normal form of $A$ (see the proof
of Theorem~\ref{herbrand}).

\begin{lem}
\label{skolemization}
If $A$ is prenex and $\modelsg A$, then $\modelsg A^H.$
\end{lem}
 
\begin{proof}
Follows from the usual laws of quantification.
\end{proof}

Our next main result will be Herbrand's theorem for $\G_V$ for $V$
uncountable. (By Theorem~\ref{bendixon}, $V$ contains a dense linear
subordering.) Let $A$ be a formula. The {\em Herbrand universe} $U(A)$
of $A$ is the set of all variable-free terms which can be constructed
from the set of function symbols occurring in $A$. To prevent $U(A)$
from being finite or empty we add a constant and a function symbol of
positive arity if no such symbols appear in~$A$.  The {\em Herbrand
base} $B(A)$ is the set of atoms constructed from the predicate
symbols in $A$ and the terms of the Herbrand universe.  In the next
theorem we will consider the Herbrand universe of a formula $\exists
\tup{x}\, A(\tup{x})$.  We fix a non-repetitive enumeration $A_1$,
$A_2$, \dots of $B(A)$, and let $X_\ell = \{\bot, A_1, \ldots,
A_\ell, \top\}$ (we may take $\top$ to be a formula which is always
$=1$).  $A(\tup{t})$ is an \emph{$\ell$-instance} of $A(\tup{x})$ if
the atomic subformulas of $A(\tup{t})$ are in $X_\ell$.

\begin{defn} 
An $\ell$-\emph{constraint} is a non-strict linear ordering $\preceq$
of $X_\ell$ s.t. $\bot$ is minimal and $\top$ is maximal.  An
interpretation~$\I$ {\em fulfils} the constraint $\preceq$ provided
for all $B, C \in X_\ell$, $B \preceq C$ iff $\I(B) \le \I(C)$.  We
say that the constraint $\preceq'$ on $X_{\ell+1}$ \emph{extends}
$\preceq$ if for all $B, C \in X_\ell$, $B \preceq C$ iff $B \preceq'
C$.
\end{defn}

\begin{prop}\label{ellI}
(a) Every $\I$ which fulfills $\preceq'$ also fulfills $\preceq$. (b)
if $\I$, $\I'$ fulfill $\preceq$, then there is a truth value
injection~$f$ on $X_\ell$, and $f(\I(A(\tup{t}))) = \I'(A(\tup{t}))$
for all $\ell$-instances $A(\tup{t})$ of $A(x)$; in particular,
$\I(A(\tup{t})) = 1$ iff $\I'(A(\tup{t})) = 1$.
\end{prop}

\begin{proof}
(a) Obvious. (b) Follows from Proposition~\ref{tviso}.
\end{proof} 

\begin{lem} \label{univH}
Let $A$ be a quantifier-free formula. If $\modelsg \exists \tup{x}\,
A(\tup{x})$ then there are tuples $\tup{t}_1, \dots \overline{t}_n$
of terms in $U(A)$, such that $\modelsg \bigvee_{i=1}^{n}
A(\tup{t}_i)$.
\end{lem}

\begin{proof}
We construct a ``semantic tree'' $\T$; i.e., a systematic
representation of all possible order types of interpretations of the
atoms $A_i$ in the Herbrand base.  $\T$ is a rooted tree whose nodes
appear at levels.  Each node at level~$\ell$ is labelled with an
$\ell$-constraint.

$\T$ is constructed in levels as follows: At level~0, the root of~$\T$
is labelled with the constraint $\bot < \top$.  Let $\nu$ be a node
added at level $\ell$ with label $\preceq$, and let $T_\ell$ be the
set of terms occurring in $X_\ell$.  Let (*) be: There is an
interpretation $\I$ that fulfils~$\preceq$ so that for some
$\ell$-instance $A(\tup{t})$, $\I(A(\tup{t}))=1$. If (*) obtains, $\nu$
is a leaf node of~$\T$, and no successor nodes are added at level
$\ell +1$.  Note that by Proposition~\ref{ellI}, any two
interpretations which fulfill $\preceq$ make the same $\ell$-instances
of $A(\tup{t})$ true; hence $\nu$ is a leaf node if and only if there
is an $\ell$-instance $A(\tup{t})$ s.t. $\I(A(\tup{t}))=1$ for
\emph{all} interpretations $\I$ that fulfil~$\preceq$.

If (*) does not obtain, for each $(\ell + 1)$-constraint $\preceq'$
extending $\preceq$ we add a successor node $\nu'$ labelled with
$\preceq'$ to $\nu$ at level $\ell+1$.

We now have two cases:

(1) $\T$ is finite. Let $\nu_1, \ldots, \nu_m$ be the leaf nodes
    of~$\T$ of levels $\ell_1$, \ldots, $\ell_m$, each labelled with a
    constraint $\preceq_1$, \dots, $\preceq_m$.  By (*), there are
    $\ell_i$-instances $A(\tup{t}_1)$, \dots, $A(\tup{t}_m)$ so that
    $\I(A(\tup{t}_i)) = 1$ for any $\I$ which fulfills $\preceq_i$. It
    is easy to see that every interpretation fulfills at least one of
    the $\preceq_i$.  Hence, for all $\I$, $\I(A(\tup{t}_1) \lor
    \ldots \lor A(\tup{t}_m)) = 1$, and so $\modelsg \bigvee_{i=1}^m
    A(\tup{t}_i)$.

(2) $\T$ is infinite.  By König's lemma, $\T$ has an infinite branch
    with nodes $\nu_0$, $\nu_1$, $\nu_2$, \dots where $\nu_\ell$ is
    labelled by $\preceq_\ell$ and is of level $\ell$.  Each
    $\preceq_{\ell+1}$ extends $\preceq_\ell$, hence we can form
    $\preceq = \bigcup_\ell \preceq_\ell$.  Let $V' \subseteq V$ be a
    non-trivial densely ordered subset of $V$, let $V' \ni c < 1$, and
    let $V'' = V' \cap [0, c)$.  $V''$ is clearly also densely
    ordered. Now let $V_c$ be $V'' \cup \{0, 1\}$, and let $h: B(A(x))
    \cup \{\bot, \top\} \to V_c$ be an injection which is so that, for
    all $A_i, A_j \in B(A(x))$, $h(A_i) \le h(A_j)$ iff $A_i \preceq
    A_j$, $h(\bot) = 0$ and $h(\top) = 1$. We define an interpretation
    $\I = \langle U(A(x)), \s\rangle$ by: $\s(f)(t_1, \ldots, t_n) =
    f(t_1, \ldots, t_n)$ for all $n$-ary function symbols~$f$ and
    $\s(P)(t_1, \ldots, t_n) = h(P(t_1, \ldots, t_n))$ for all $n$-ary
    predicate symbols~$P$ (clearly then, $\I(A_i) = h(A_i)$). By
    definition, $\I$ $\ell$-fulfills $\preceq_\ell$ for all $\ell$. By
    (*), $\I(A(\tup{t})) < 1$ for all $\ell$-instances $A(\tup{t})$ of
    $A(x)$, and by the definition of $V_c$, $\I(A(\tup{t})) < c$.
    Since every $A(\tup{t})$ with $\tup{t} \in U(A(x))$ is an
    $\ell$-instance of $A(x)$ for some $\ell$, we have $\I(\exists x\,
    A(\tup{x})) \le c < 1$.This contradicts the assumption that
    $\modelsg \exists \tup{x}\, A(\tup{x})$.
\end{proof}

The following lemma establishes sufficient conditions for a logic to
allow {\em reverse Skolemization}. By this we mean the re-introduction
of quantifiers in Herbrand expansions.
Here, by a {\em logic}~$\lo$ we mean a set of formulas that is closed
under modus ponens, generalization and substitutions (of both formulas
and terms).  We call a formula $A$ {\em valid} in $\lo$, $\modelsl A$,
if $A \in \lo$.  The following three results follow from
\cite{BCF01LPAR} together with Lemma~\ref{univH}:

\begin{lem}\label{propreskolem}
Let $\lo$ be a logic satisfying the following properties:

(1) $\modelsl A \lor B  \IMPLIES 
        \modelsl B \lor A$

(2) $\modelsl (A \lor B) \lor C \IMPLIES 
         \modelsl A \lor (B \lor C)$ 

(3) $\modelsl A \lor B \lor B \IMPLIES 
       \modelsl A \lor B$ 

(4) $\modelsl A(y) \IMPLIES \modelsl \forall x\, A(x)$ 

(5) $\modelsl A(t) \IMPLIES \modelsl \exists x\,  A(x) $

(6) $\modelsl  \forall x (A(x) \lor B) \IMPLIES 
\modelsl \forall x\, A(x) \lor B$

(7) $\modelsl  \exists x (A(x) \lor B) \IMPLIES 
\modelsl \exists x\, A(x) \lor B$.

\noindent ($x$ is not free in $B$.)  Let $\exists \tup{x}\,
A^F(\tup{x})$ be the Herbrand form of the prenex formula $\tup{\Q}_i
A(\tup{y}_i)$, and let $\tup{t}_1, \ldots, \tup{t}_m$ be tuples of terms in
$U(A^F(\tup{x}))$. If $\modelsl \bigvee_{i=1}^{m} A^F(\tup{t_i})$, then
$\modelsl \tup{\Q}\tup{y}\, A(\tup{y}).$
\end{lem}

\begin{cor}
\label{skolreskol}  
If $\modelsg \exists \tup{x}\, A^F(\tup{x})$, then $\modelsg \tup{\Q}
\tup{y} A(\tup{y}).$
\end{cor}

\begin{thm} \label{HT}
Let $A \equiv \tup{\Q} \tup{y} B(\tup{y})$ be prenex. $\modelsg
\tup{\Q} \tup{y} B(\tup{y})$ iff there are tuples
$\tup{t}_1, \ldots \tup{t}_m$ of terms in $U(A^H(\tup{x}))$, such that
$\modelsg \bigvee_{i=1}^{m} B^F(\tup{t}_i).$
\end{thm}

\begin{rem} 
An alternative proof of Herbrand's theorem can be obtained using the
analytic calculus {\it HIF} (``Hypersequent calculus for
Intuitionistic Fuzzy logic'') \cite{BaazZach00CSL}.
\end{rem}

\begin{thm}
The prenex fragment of a Gödel logic based on a truth value set~$V$
which is either finite or uncountable infinite is axiomatizable. An
axiomatization is given by the standard axioms and rules for~$\LC$
extended by conditions (4)--(7) of Lemma~\ref{propreskolem} written as
rules. For the $m$-valued case add the characteristic axiom for $\G_m$,
\(
G_m \equiv \bigvee_{i=1}^m\bigvee_{j=i+1}^{m+1} ((A_i \impl A_j) \land 
(A_j \impl A_i)).
\)
%Logics of truth value sets with the same cardinality are the same; and
%for quantifier-free formulas, \LC{} and $\LC + G_n$ are complete for
%$\G_V$ if $V$ is infinite or $n$-valued, respectively 
\end{thm}

\begin{proof}  
Completeness: Let $\tup{\Q} \tup{y}_i A(\tup{y})$ be a prenex formula
valid in $\G_V$. Herbrand's theorem holds for $\G_V$ (for $V$
infinite, this is Theorem~\ref{HT}; for $V$ finite it follows from
results in \cite{BaazFermZach94IPC}), and so a Herbrand disjunction
$\bigvee_{i=1}^{n} A^F(\tup{t}_i)$ is provable in \LC{} or $\LC + G_m$
\cite[Chapter 10.1]{gottwald}.  $\tup{\Q} \tup{y}_i A(\tup{y})$ is
provable by Lemma~\ref{propreskolem}.

Soundness: $\G_V$ satisfies the conditions of Lemma~\ref{propreskolem}
(in particular, note that $\forall x(A(x) \lor B) \to (\forall x\,
A(x) \lor B)$ with $x$ not free in $B$ is valid in all G\"odel
logics).
\end{proof}

%%%%%%%%%%%%%%%%%%%%%%%%%%%%%%%%%%%%%%%%%%%%%%%%%%%%%%%%%%%%%%%%%%%%%%%%%
%
% NON AXIOMATIZABILITY RESULTS
%
%%%%%%%%%%%%%%%%%%%%%%%%%%%%%%%%%%%%%%%%%%%%%%%%%%%%%%%%%%%%%%%%%%%%%%%%%

\section{Nonaxiomatizability results}

In this section we show that the prenex fragments of first-order Gödel
logics where the set of truth values does not contain a dense subset
are not axiomatizable. We establish the result first for the entire
set of valid formulas by reducing the classical validity of a formula
in all finite models to the validity of a formula in Gödel logic (the
set of these formulas is not r.e.\ by Trakhtenbrot's Theorem).  We
then strengthen the result by showing that the image of the
translation from the prenex fragment of classical logic to Gödel logic
is equivalent to a prenex formula.

\begin{thm}\label{nonax}
If $V$ is countably infinite, then $\G_V$ is
not axiomatizable.
\end{thm}

\begin{proof}
By Theorem~\ref{bendixon}, $V$ is countably infinite iff it is
infinite and does not contain a non-trivial densely ordered subset.
We show that for every sentence $A$ there is a sentence $A^g$
s.t. $A^g$ is valid in $\G_V$ iff $A$ is true in every finite
(classical) first-order structure.

We define $A^g$ as follows: Let $P$ be a unary and $L$ be a binary
predicate symbol not occurring in~$A$ and let $Q_1$, \dots, $Q_n$ be
all the predicate symbols in~$A$.  We use the abbreviations $x \in y
\equiv \neg\neg L(x, y)$ and $x \prec y \equiv (P(y) \impl P(x)) \impl
P(y)$.  Note that for any interpretation~\I, $\I(x \in y)$ is either
$0$ or $1$, and as long as $\I(P(x)) < 1$ for all $x$ (in particular,
if $\I(\exists z\, P(z)) < 1$), we have $\I(x \prec y) = 1$ iff
$\I(P(x)) < \I(P(y))$.  Let $A^g \equiv$
\[
\left\{
\begin{array}{l}
S \land {} c_1 \in 0 \land {} c_2 \in 0 \land 
c_2 \prec c_1 \land {}\\
\quad \forall i\bigl[\forall x, y \forall j \forall k \exists z\, \succ
\lor \forall x\neg(x \in i)\bigr]
\end{array}\right\}
\impl (A' \lor \exists u\, P(u))
\]
where $S$ is the conjunction of the standard axioms for $0$, successor
and $\le$, with double negations in front of atomic formulas,\vskip-15pt
$$\succ\equiv \begin{array}{l}
(j \le i \land x \in j \land k\le i \land y \in k \land x \prec y) \impl {}\\
\qquad \impl (z \in s(i) \land x \prec z \land z \prec y)
\end{array}$$
and
$A'$ is $A$ where every atomic formula is replaced by its double
negation, and all quantifiers are relativized to the predicate
$R(i)\equiv \exists x(x \in i)$.

Intuitively, $L$ is a predicate that divides a subset of the domain
into levels, and $x \in i$ means that $x$ is an element of level~$i$.
$P$ orders the elements of the domain which fall into one of the
levels in a subordering of the truth values.  The idea is that for any
two elements in a level $\le i$ there is an element in level $i + 1$
which lies strictly between those two elements in the ordering given
by~$\prec$.  If this condition cannot be satisfied, the levels above
$i$ are empty.  Clearly, this condition can be satisfied in an
interpretation~\I{} only for finitely many levels if $V$ does not
contain a dense subset, since if more than finitely many levels are
non-empty, then $\bigcup_{i} \{\I(P(d)) : \I \models d \in i\}$ gives
a dense subset.  By relativizing the quantifiers in $A$ to the indices
of non-empty levels, we in effect relativize to a finite subset of the
domain.  We make this more precise:

Suppose $A$ is classically false in some finite structure~\I.
W.l.o.g. we may assume that the domain of this structure is the
naturals $0$, \dots, $n$.  We extend $\I$ to a $\G_V$-interpretation
$\I^g$ with domain $\N$ as follows: Since $V$ contains infinitely many
values, we can choose $c_1$, $c_2$, $L$ and $P$ so that $\exists x(x
\in i)$ is true for $i = 0$, \dots, $n$ and false otherwise, and so
that $\sup \Distr_{\I^g} P(x) < 1$.  The number-theoretic symbols receive
their natural interpretation.  The antecedent of $A^g$ clearly receives
the value~1, and the consequent receives $\sup \Distr{\I^g} P(x) < 1$,
so $\I^g \nmodels A^g$.

Now suppose that $\I \nmodels A^g$.  Then $\I(\exists x\, P(x)) < 1$
and so $\sup \Distr_\I P(x) < 1$.  In this case, $\I(x \prec y) = 1$
iff $\I(P(x)) < \I(P(y))$, so $\prec$ defines a strict order on the
domain of $\I$.  It is easily seen that in order for the value of the
antecedent of $A^g$ under \I{} to be greater than that of the
consequent, it must be $= 1$ (the values of all subformulas are either
$\le \sup \Distr_\I P(x)$ or $= 1$).  For this to happen, of course,
what the antecedent is intended to express must actually be true
in~$\I$, i.e., that $x \in i$ defines a series of disjoint levels and
that for any $i$, either level $i + 1$ is empty or for all $x$, $y$
s.t. $x \in j$, $y \in k$ with $j, k \le i$ and $x \prec y$ there is
a~$z$ with $x \prec z \prec y$ and $z \in i+1$.  To see this, consider the 
relevant part of the antecedent, $B = \forall i\bigl[
\forall x, y\forall j\forall k\exists z\,\succ 
\lor \forall x\neg(x \in i)\bigr]$.
If $\I(B) = 1$, then for all $i$, either
\(
\I(\forall x, y\forall j\forall k\exists z\,\succ) = 1
\)
or $\I(\forall x\neg(x \in i)) = 1$. In the first case, we have
\(
\I(\exists z\,\succ) = 1
\)
for all $x$, $y$, $j$, and $k$. Now suppose that
for all $z$, $\I(\succ) < 1$, yet $\I(\exists z\, \succ) = 1$.
Then for at least some $z$ the value of that formula would have to be
$> \sup \Distr_\I P(z)$, which is impossible.  Thus, for every $x$,
$y$, $j$, $k$, there is a $z$ such that $\I(\succ) = 1$.  But this
means that for all $x$, $y$ s.t. $x \in j$, $y \in k$ with $j, k \le
i$ and $x \prec y$ there is a~$z$ with $x \prec z \prec y$ and $z \in
i+1$.

In the second case, where $\I(\forall x\neg(x \in i)) = 1$, we have
that $\I(\neg (x \in i)) = 1$ for all $x$, hence $\I(x \in i) = 0$ and
level~$i$ is empty.

Since $V$ contains no dense subset, from some finite level $i$ onward,
the levels must be empty. Of course, $i > 0$ since $c_1 \in 0$.  Thus, $A$
is false in the classical interpretation $\I^c$ obtained from
$\I$ by restricting $\I$ to the domain $\{0, \ldots, i - 1\}$ and
$\I^c(Q) = \I(\neg\neg Q)$ for atomic $Q$.
\end{proof}

This shows that no infinite-valued Gödel logic whose set of truth
values does not contain a dense subset is axiomatizable.  We
strengthen this result to show that the prenex fragments are likewise
not axiomatizable.  This is done by showing that if $A$ is prenex,
then there is a formula $A^G$ which is also prenex and which is valid
in $\G_V$ iff $A^g$ is.  Since not all quantifier shifting rules are
generally valid, we have to prove that in this particular instance there is a
prenex formula which is valid in $\Gd$ iff $A^g$ is.

\begin{thm}
If $V$ is countably infinite, the prenex fragment of $\G_V$ is not
axiomatizable.
\end{thm}

\begin{proof}
By the proof of Theorem~\ref{nonax}, a formula $A$ is true in all
finite models iff $\G_V \models A^g$.  $A^g$ is of the form $B \impl
(A' \lor \exists u\, P(u))$.  We show that $A^g$ is equivalent in
$\G_V$ to a prenex formula.  

Call a formula $A$ in which every atomic formula occurs negated a
\emph{classical} formula.  It is easy to see that for any $\I$ and
$A(x)$ with $\I(A(d)) \in \{0, 1\}$ for all $d$, $\I(\forall x\, A(x)
\to B) = \I(\exists x\, A(x) \to B)$ and $\I(B \to \exists x\, A(x)) =
\I(\exists x(B \to A(x))$. Hence, any classical formula is equivalent
to a prenex formula; let $A_0$ be a prenex form of $A'$.  Since
all quantifier shifts for conjunctions are valid, the antecedent~$B$
of $A^g$ is equivalent to a prenex formula $\Q_1x_1 \ldots \Q_n x_n
B_0(x_1, \ldots, x_n)$.  Hence, $A^g$ is equivalent to
$\tup{\Q}\tup{x} B_0(\tup{x}) \to (A_0 \lor \exists u\, P(u))$.

Let $\Q_i'$ be $\exists$ if $\Q_i$ is $\forall$, and $\forall$ if
$\Q_i$ is $\exists$, let $C \equiv A_0 \lor \exists u\, P(u)$, and $v
= \I(\exists u\, P(u)))$.  We show that $\tup{\Q}\tup{x}\,
B_0(\tup{x}) \to C$ is equivalent to $\tup{\Q}'\tup{x}(B_0(\tup{x})
\to C)$ by induction on~$n$.  Let $\tup{\Q}\tup{x}B_0 \equiv \Q_1
x_1\ldots \Q_{i}x_{i} B_1(d_1, \ldots, d_{i-1}, x_i)$.  Since
quantifier shifts for $\exists$ in the antecent of a conditional are
valid, we only have to consider the case $\Q_i = \forall$. Suppose
$\I(\forall x_i\, B_1(\tup{d}, x_i) \impl C) \neq \I(\exists
x_i(B_1(\tup{d}, x_i) \impl C)$. This can only happen if $\I(\forall
x_i\, B_1(\tup{d}, x_i)) = \I(C) < 1$ but $\I(B_1(\tup{d}, c)) > \I(C)
\ge v$ for all $c$.  However, it is easy to see by inspecting $B$ that
$\I(B_1(\tup{d}, c))$ is either $=1$ or $\le v$.

Now we show that $\I(B_0(\tup{d}) \impl (A_0 \lor \exists u\, P(u))) =
\I(\exists u(B_0(\tup{d}) \impl (A_0 \lor P(u))))$.  If $\I(A_0) = 1$,
then both sides equal $= 1$.  If $\I(A_0) = 0$, then $\I(A_0 \lor
\exists u\, P(u)) = v$.  The only case where the two sides might
differ is if $\I(B_0(\tup{d})) = v$ but $\I(A_0 \lor P(c)) = \I(P(c))
< v$ for all $c$.  But inspection of $B_0$ shows that $\I(B_0(\tup{t}))
= 1$ or $= \I(P(e))$ for some $e \in \tup{d}$ (the only subformulas of
$B_0(\tup{d})$ which do not appear negated are of the form $e' \prec
e$).  Hence, if $\I(B_0(\tup{d})) = v$, then for some~$e$, $\I(P(e)) =
v$.

Last we consider the quantifiers in $A_0 \equiv \tup{\Q} \tup{y}\, A_1$.
Since $A_0$ is classical, $\I(B_0(\tup{d}) \impl (A_0 \lor P(c))) =
\I(\tup{\Q}\tup{y}(B_0(\tup{d}) \impl (A_1 \lor P(c))))$ for all
$\tup{d}$, $c$.  To see this, first note that shifting quantifiers
across $\lor$, and shifting universal quantifiers out of the
consequent of a conditional is always possible. Hence it suffices to
consider the case of $\exists$.  $\I(\exists y\, A_2)$ is
either $= 0$ or $= 1$.  In the former case, both sides equal
$\I(B_0(\tup{d}) \impl P(d))$, in the latter, both sides equal~$1$.
\end{proof}
In summary, we obtain the following characterization of axiomatizability of
prenex fragments of Gödel logics:
\begin{thm}
The prenex fragment of $\G_V$ is axiomatizable if and only if $V$ is
finite or uncountable.
\end{thm}

%%%%%%%%%%%%%%%%%%%%%%%%%%%%%%%%%%%%%%%%%%%%%%%%%%%%%%%%%%%%%%%%%%%%%%%%%
%
% CONSEQUENCES
%
%%%%%%%%%%%%%%%%%%%%%%%%%%%%%%%%%%%%%%%%%%%%%%%%%%%%%%%%%%%%%%%%%%%%%%%%%

\section{Conclusion}

Our characterization relates in an interesting way to compactness
results of entailment relations of Gödel logics as given
in~\cite{BaazZach:98}: Exactly those Gödel logics have an
axiomatizable prenex fragment which also have a compact propositional
logic.

For full first order Gödel logics the situation is quite similar in
the sense that the truth value set must be finite or uncountable to
allow axiomatization, but in addition it is necessary that 0 either be
in the perfect kernel of the truth value set or be isolated. Two
different logics correspond to these conditions, which have the same
prenex fragment. Consequently there are Gödel logics where the prenex
fragment is axiomatizable, but the full logic does not allow a
recursive axiomatization. These are the logics of truth value sets
which contain an uncountable subset, but 0~is neither in the perfect
kernel nor isolated. These results have been obtained
in~\cite{PreiDiss} and will be reported in a forthcoming article by
the authors.

\end{document}